\documentclass[11pt]{article}
\usepackage{epsfig}
\usepackage{t1enc}
    \usepackage[latin1]{inputenc}
    \usepackage[english]{babel}
        \usepackage{latexsym}
\usepackage{amssymb}
\usepackage{euscript}

\newcommand{\bs}{\textbf{S}}

\def\begg{\begin{equation}}
\def\endd{\end{equation}}

\newcommand{\bp}{{\bf P}}
\newcommand{\ep}{{\epsilon}}
\begin{document}
\setlength{\arraycolsep}{.136889em}
\renewcommand{\theequation}{\thesection.\arabic{equation}}
\newtheorem{thm}{Theorem}[section]
\newtheorem{propo}{Proposition}[section]
\newtheorem{lemma}{Lemma}[section]
\newtheorem{corollary}{Corollary}[section]
\newtheorem{remark}{Remark}[section]
\centerline{\Large\bf HOW TALL CAN BE THE EXCURSIONS}

\smallskip

\centerline{\Large\bf OF A RANDOM WALK ON A SPIDER?}

\bigskip\bigskip
\bigskip\bigskip
\centerline{\it Dedicated to  the memory of Marc Yor. }

\bigskip\bigskip
\bigskip\bigskip
\renewcommand{\thefootnote}{1}

\noindent
{\textbf{Ant\'{o}nia
F\"{o}ldes}\footnote{Research supported by a PSC CUNY Grant, No.
68030-0043.}}
\newline
Department of Mathematics, College of Staten
Island, CUNY, 2800 Victory Blvd., Staten Island, New York 10314,
U.S.A.  E-mail address: antonia.foldes@csi.cuny.edu

\bigskip
\noindent
\renewcommand{\thefootnote}{2}
\noindent
{\textbf{Pal R\'ev\'esz}\footnote{Research supported by the
Hungarian National Foundation for Scientif\/ic Research, Grant No.
K 108615}}
\newline
\noindent
Institut f\"ur
Statistik und Wahrscheinlichkeitstheorie, Technische Universit\"at
Wien, Wiedner Hauptstrasse 8-10/107 A-1040 Vienna, Austria.
E-mail address: reveszp@renyi.hu
\bigskip

\medskip
\noindent{\bf Abstract}\newline
We consider a simple symmetric random walk on a spider, that is a collection of half lines (we call them legs) joined at the origin. Our main question is the following: if the walker makes $n$ steps how high can he go up on all legs. This problem is discussed in two different situations; when the number of legs are increasing, as $n$ goes to infinity and when it is fixed.


\medskip

\vspace{.1cm}

\medskip

\section{Introduction and main results}
\renewcommand{\thesection}{\arabic{section}} \setcounter{equation}{0}
\setcounter{thm}{0} \setcounter{lemma}{0}

 In 1978 Walsh \cite{W} introduced  a Brownian motion which lives on $N$-semiaxis joined at the origin. This motion loosely speaking
 performs regular Brownian motion on each semiaxis,  and when it arrives to the origin it continues its motion
on any of the $N$ semiaxis with  equal probability. This "definition" can be made precise  with excursion theory or with the computation of the semigroup of this motion, see e.g. Barlow, Pitman and Yor  \cite{BPY89}. This motion is called   now  Walsh's Brownian motion, or Brownian spider. In \cite{PPL} the exit time from specific sets is investigated, and a generalized arc-sine law is introduced on the time spent on a specific semiaxis. This question is further investigated in the elegant paper of  Vakeroudis  and Yor (2012) \cite{VY}.
In 2013 the second author of this paper introduced in his book \cite{R13} the discrete version of the above motion, and called it a random walk on a spider. His main interest however was to consider this spiderwalk where the number of legs of the spider goes to infinity. In this paper we consider the spiderwalk in both situations, with  finite  and increasing numbers of legs.
\smallskip

 We  start with  the definitions, lifted form the book:
Let $\textbf{SP(N)}=(V_N,E_N),$ where

$$V_N=\left\{v_N(r,j)=r\exp\left(\frac{2\pi ij}{N}\right) \quad r=0,1,...,\quad j=1,...,N\right\}\quad \rm{and} \,\,i=\sqrt{-1}.$$
\noindent
is the set of vertices of \textbf{SP(N)} and

$$E_N=\{e_N(r,j)=(v_N(r,j),v_N(r+1,j)  \quad r=0,1,...,\,\,j=1,...,N \}$$
\noindent
is the set of edges of \textbf{SP(N)}. We will call \textbf{SP(N)}  a spider with $N$ legs. The vertex

$$v_N(0)=v_N(0,1)=v_N(0,2)=...=v_N(0,N)$$
\noindent
is called the body of  the spider, while   $\{v_N(1,j), v_N(2,j)...\}$ is the $j$-th leg of the spider.

On \textbf{SP(N)}  we consider a random walk $\{\bs_n\}_{n=0}^{\infty}$ starting from the body of spider $\bs_0=v_N(0),$ with
 the following transition probabilities:
$${\bf P} (\bs_{n+1}=v_N(1,j)|\bs_n=v_N(0))=\frac{1}{N}\quad j=1,...,N $$
and for $r=1,...,\quad j=1,...,N$
$${\bf P} (\bs_{n+1}=v_N(r+1,j)|\bs_n=v_N(r,j))={\bf P} (\bs_{n+1}=v_N(r-1,j)|\bs_n=v_N(r,j))=\frac{1}{2}.$$
Let
\begin{eqnarray*}
\xi(v_N(r,j),n):&=&\#\{k: k\leq n,\, \bs_k=v_N(r,j)\} \\
\zeta(n):&=&\#\{k: k\leq n,\, \bs_k=v_N(0)\}=\xi(v_N(0),n)
\end{eqnarray*}
and define the events

\begin{eqnarray*}
M(n,R):&=&\{  \min_{1\leq j \leq N} \xi(v_N(R,j),n)\geq 1 \} \\
A(n,R,k):&=&\{  \min_{1\leq j \leq N} \xi(v_N(R,j),n)\geq k \}.
\end{eqnarray*}
Observe that the meaning of the event  $M(n,R)$ is that in $n$ steps the walker climbs up to at least $R$  on each legs.
The special case $M(n,1)$ means that in $n$ steps each leg is visited at least once. $A(n,R,k)$ means that in $n$ steps the walker visits each legs at height $R$ at least $k$ times.
\newpage
We recall the main result from the book;
\smallskip

\noindent
{\bf Theorem A:} {\it For the}  \textbf{SP(N)}

\begg
\lim_{N\to \infty} {\bf P}( M(N\log N)^2,1))= \left(\frac{2}{\pi}\right)^{1/2} \int_1^\infty e^{-u^2/2}\, du= {\bf P}(|Z|>1).
\endd
{\it where} $Z$  {\it is a standard normal random  variable.}

\smallskip
In this paper we ask what can we say about $M(n,L).$ Our main result is
\begin{thm}  For any integer $L \leq \frac{N}{\log N}$  we have for the {\rm \textbf{SP(N)}}
\begg
\lim_{N\to \infty} {\bf P}( M(cLN\log N)^2,L))={\bf P}\left(|Z|>\frac{1}{c}\right):=p(c).
\endd
\end{thm}

\smallskip
\noindent
To formulate in  words, the theorem above gives the limiting probability of the event that as $N\to \infty,$    in $(cLN\log N)^2$
 steps the walker arrives at least at height  $L$ on each of the N legs at least once.

\noindent
The next two results are natural companions of the above one;
\begin{thm}  For any integer $L \leq \frac{N}{\log N}, $ {\it  and any sequence} $f(N)\uparrow \infty,$  we have for the {\rm\textbf{SP(N)}}
\begg
\lim_{N\to \infty} {\bf P}( M(f(N)L N\log N)^2,L))=1.
\endd
\end{thm}

 \begin{thm}  For any integer $L \leq \frac{N}{\log N}, $  and any sequence $f(N)\downarrow 0, $ we have for the {\rm\textbf{SP(N)}}
\begg
\lim_{N\to \infty} {\bf P}( M(f(N)L N\log N)^2,L))=0.
\endd
\end{thm}

\noindent
Furthermore we have
\begin{thm}  For any integer $L \leq \frac{N}{\log N}$ and any fixed integer $k \geq 1,$ we have for the {\rm \textbf{SP(N)}}
\begg
\lim_{N\to \infty} {\bf P}( A(cLN\log N)^2,L,k))={\bf P}\left(|Z|>\frac{1}{c}\right):=p(c).
\endd
\end{thm}

\begin{thm}  For any integer $L \leq \frac{N}{\log N},$   and any fixed integer $k \geq 1,$   {\it  and any sequence} $f(N)\uparrow \infty,$  we have for the {\rm \textbf{SP(N)}}
\begg
\lim_{N\to \infty} {\bf P}( A(f(N)L N\log N)^2,L,k))=1.
\endd
\end{thm}

\noindent
In the above theorems the $L \leq \frac{N}{\log N} $  condition is a technical one, which should be eliminated. So we ask the following
\newpage
\noindent
{\bf Question 1:} {\it Determine for each $0\leq p\leq 1,$ the function $g(N,L,p)$  such that for the {\rm \textbf{SP(N)}}
$$\lim_{N\to \infty} {\bf P}( M(g(N,L,p),L))=p$$
should hold.}
\bigskip

\noindent
{\bf Question 2:} {\it Determine for each $0\leq p\leq 1,$ the function $g^*(N,L,p)$  such that for the {\rm \textbf{SP(N)}}
$$\lim_{N\to \infty} {\bf P}( A(g^*(N,L,p),L,k))=p$$
should hold.}

\noindent

\bigskip

In the second part of this paper we will consider the spider with $K$ legs where $K$ is a fixed integer.  Of course the $K=2$ case is the simple symmetric walk on the line.
One of the natural questions to ask is how high  does the walker go up on the legs of this spider. Using the the definitions above for $K$ legs,
denote

\begin{eqnarray*}
M_K(n,j) &=& \max\{r; \xi(v_K(r,j),n) \geq 1\}, \\
M_K(n) &=& \max_{1\leq j\leq K} M_K(n,j) \\
\end{eqnarray*}
Clearly $M_K(n,j)$ is the highest point on leg $j, $ where the walker arrives in $n$ steps.
The LIL for the simple random walk clearly implies that

$$\limsup_{n\to\infty} \frac{M_K(n)}{\sqrt{2n\log\log n}}=1. \quad {\rm a.s.}$$
as $M_K(n)$ is the same for $\bf{SP(K)}$ and for the simple symmetric walk. Similarly from the other LIL (see Lemma E in the next section) we also have that

$$\liminf_{n\to\infty} \left(\frac{\log\log n}{n}\right)^{1/2}M_K(n)=\frac{\pi}{\sqrt{8}}. \quad {\rm a.s.}$$

  However it is a much more interesting question to ask  the maximal height which can be achieved on all legs simultaneously.  To be  more precise we ask what can we  say about
 $\min_{1 \leq j \leq K} M_K(n,j).$  We will prove that

 \begin{thm} For the $\bf{SP(K)}$ we have
\begg
\limsup_{n\to\infty} \frac {\min_{1 \leq j\leq K} M_K(n,j).}{\sqrt{2n\log\log n}}=\frac{1}{2K-1} \quad a.s.
\endd
 \end{thm}

\section{Preliminary  Results}
\renewcommand{\thesection}{\arabic{section}} \setcounter{equation}{0}
\setcounter{thm}{0} \setcounter{lemma}{0}

\noindent
We will need the famous Erd\H{o}s- R\'enyi  \cite{ER61} coupon collector  theorem:

\smallskip
\noindent
{\bf Theorem B:} {\it Suppose that there are $N$ urns given, and that $N\log N+(m-1)N\log\log N+Nx$ balls are placed in these urns one after the other independently. Then for every real $x$ the probability that each urn
 will contain at least m balls converges to}

\begg
 \exp\left( -\frac{1}{(m-1)!}\exp(-x)\right),
\endd
as
$N\to\infty.$

 \noindent
It is worthwhile to spell out the most important special case $m=1,$ as follows;

\smallskip
\noindent
{\bf Theorem C:} {\it Suppose that there are $N$ urns given, and that $N\log N+Nx$ balls are placed in these urns one after the other independently. Then for every real $x$ the probability that each urn
will contain at least one ball converges to}

\begg
 \exp\left( -\exp(-x)\right),
 \endd
as $N\to\infty.$

\noindent
We will also need Hoeffding \cite{H61} inequality;

\smallskip
\noindent
{\bf LEMMA D}: { \it Let} $a_i\leq X_i \leq b_i \quad(i=1,2,...k)$ { \it be independent random variables and } $S_k=\sum _{i=1}^k X_i.$ { \it Then for every} $x>0$

\begg
{\bf P}(|S_k - {\bf E} S_k|\geq k x)\leq 2\exp\left(-\frac{2k^2x^2}{\sum_{i=1}^k(b_i-a_i)^2}\right). \label{ho1}
\endd
We will use the above inequality  in the following special case:

Let $X_1, X_2,...X_j\,$ i.i.d. Bernoulli random variables, then for $j\leq k$
\begin{equation}
{\bf P}(|S_j - {\bf E} S_j|\geq k x)\leq 2\exp\left(-2kx^2\right). \label{ho2}
\end{equation}
To see this, enough to observe that for $j\leq k$ we might take $X_{j+1}=X_{j+2}=...=X_{k}=0$, then

\noindent
 $\sum_{i=1}^n(b_i-a_i)^2=j.$

\noindent
Recall the famous other LIL of Chung  for the simple symmetric random walk see e.g. in \cite{R81} (page 41)

\smallskip
\noindent
 {\bf Lemma E}: {\it For the  the simple symmetric random walk} $\{S_i\}_{i=0}^\infty$ { \it let} $$M(n)=\max_{0\leq i\leq n} |S_i|.$$ {\it Then we have}
 \begg
 \liminf_{n \to \infty} \left(\frac{\log\log n}{n} \right)^{1/2} M(n)=\frac{\pi}{\sqrt{8}} \quad {\rm a.s.} \label{oi}
 \endd

\noindent
We will use the celebrated functional law of iterated logarithm of Strassen.
Let ${\cal S}$ be the Strassen class of functions, i.e.,
${\cal S}\subset C([0,1],{\mathbb R})$ is the class of absolutely
continuous functions (with respect to the Lebesgue measure) on $[0,1]$
for which
\begin{equation}
f(0)=0\qquad {\rm and\qquad } I(f)=\int_0^1\dot{f}^2(x)dx\leq 1.
\end{equation}
\noindent
Denote the continuous versions of the random walk process $\{S(nx);\,\, 0\leq x \leq 1\}_{n=1}^{\infty}$ defined  by linear interpolation from the simple symmetric random walk $\{S_n\}_{n=0}^{\infty}.$
\medskip\noindent

\noindent
{\bf Theorem F} \cite{ST}
{\it The sequence of random functions
$$
\left\{\frac{S(xn)}{(2n\log\log n)^{1/2}};\, 0\leq x\leq
1\right\}_{n\geq 3},
$$
as $n\to\infty$, is almost surely relatively compact in the space
$C([0,1])$ and the set of its limit points is the class of
functions ${\cal S}$.}

\noindent
Recall the following well-known lemma (see e.g. in \cite{FR})

\noindent
 {\bf Lemma G}: {\it Let

 $$0=t_0\leq t_1\leq...\leq t_n=1$$
 Define $g(.)\in C[0,1]$ by the requirements:
 $$g(t_m)=f(t_m) \,\, {\rm for} \,\, m=0,1,...,n$$
 $$g(.) \,\,{\rm is\,\,linear\,\,on}\,\, [t_{m-1},t_m]\,\,{\rm for}\,\, m=1,...,n$$
then
$$\int_0^1\dot{g}^2(t)dt\leq\int_0^1\dot{f}^2(t)dt.$$
Equality holds iff } $g(t)\equiv f(t).$

\noindent
The next lemma is also well-known, an easy reference is  \cite{CCFR97}.

\smallskip
\noindent
 {\bf Lemma H}: {\it Let $\{A_k\}_{k\geq 1}$  be an arbitrary sequence of events such  that we have $\bp(A_k \,\,i.o.)=1$.
 Let $\{B_k\}_{k\geq 1}$ be another arbitrary sequence of events  that is independent of $\{A_k\}_{k\geq 1}$ an assume that
$\bp(B_k \,\,i.o.)\geq p>0.$  Then we have also} $\bp(A_k\,B_k \, \,i.o.)=1.$
\newpage
 \section{Proofs}
\renewcommand{\thesection}{\arabic{section}} \setcounter{equation}{0}
\setcounter{thm}{0} \setcounter{lemma}{0}

Let $\{S_n\}_{n=0}^{\infty}$ be a  simple symmetric  one dimensional random walk and let
\begin{eqnarray*}
\xi(0,n)&=&\#\{k: 1\leq k<n,\ S_k=0\},\\
\zeta(1,n)&=&\#\{k: 1\leq k<n,\ S_k=0,\ |S_{k+1}|=1\},\\
\zeta(L,n)&=&\#\{k:1\leq k<n,\ S_k=0\ {\rm and\,}  |S_{k+i}|\,\,\,i=1,2... \,{\rm  hits}\,\, L \,\,{\rm before \ returning \ to \,  }\, 0
\},\\
\rho(0)&=&0\, \ {\rm and\,\,} \rho(m)=\min\{k: k>\rho(m-1),\, S_k=0\} \\
\end{eqnarray*}
Then $\xi(0,\rho_m)=m.$
Finally let
$$ H(n)=\rho(\xi(0,n)+1).$$

\begin{lemma}
\begin{equation}
|\zeta(L,H(n))-\zeta(L,n)|\leq 1 \quad a.s.
\end{equation}
\begin{equation}
|\xi(0,H(n))-\xi(0,n)|\leq 1\quad a.s.
\end{equation}
\end{lemma}

\vspace{2ex}\noindent
{\bf Proof:} Trivial.

\begin{lemma}
$$\bp\left(|\zeta(L,n)-L^{-1}\xi(0,n)|\geq  4 n^{1/4}(\log n)^{3/4}\right)\leq \frac{2}{n}$$
\end{lemma}
for $n$ big enough.
\vspace{2ex}\noindent

\noindent
{\bf Proof:} Let

$$D(n)=|\zeta(L,H(n)) - L^{-1} \xi(0, H(n))|=
              |\zeta(L,\rho(\xi(0,n)+1)) -L^{-1} \xi(0,\rho(\xi(0,n)+1))|.$$
As
$$|\zeta(L,n)-L^{-1}\xi(0,n)|\leq |\zeta(L,H(n))-\xi(0,H(n))|+2,$$
we get for $n$ big  enough, that
\begin{eqnarray*}
\lefteqn{\bp(|\zeta(L,n)-L^{-1}\xi(0,n)|\geq 4 n^{1/4}(\log n)^{3/4})}\\
&\leq &\bp(|\zeta(L,H(n))-L^{-1}\xi(0,H(n))|\geq 3 n^{1/4}(\log n)^{3/4})=\\
&= &\bp(|D(n)|\geq 3\, n^{1/4}(\log n)^{3/4},\,\xi(0,n)\geq2 n^{1/2} (\log n)^{1/2})+  \\
&+ & \bp(|D(n)|\geq 3\, n^{1/4}(\log n)^{3/4},\,\xi(0,n)< 2 n^{1/2} (\log n)^{1/2})=\\
&=& I +II.
\end{eqnarray*}

We have for $n$ big enough, that
\begin{eqnarray*}
 I&\leq& \bp(\xi(0,n)\geq 2 n^{1/2} (\log n)^{1/2})\leq \bp(|Z|\geq \frac{3}{2}\,(\log n)^{1/2} )\leq \frac{1}{n}, \\
 II &\leq& \sum_{i=1} ^{2 n^{1/2} (\log n)^{1/2}} \bp(|\zeta(L,\rho(i)) -L^{-1}\xi(0,\rho(i))| >
3 n^{1/4} (\log n)^{3/4},\, \xi(0,n)=i )\\
&\leq& \sum_{i=1} ^{2 n^{1/2} (\log n)^{1/2}} \bp(|\zeta(L,\rho(i)) -L^{-1}i)| >3\, n^{1/4}(\log n)^{3/4})\\
&\leq& \sum_{i=1} ^{2 n^{1/2} (\log n)^{1/2}} 2 \exp(-9\,\log n) \\
& \leq & 4 n^{1/2} (\log n)^{1/2} \exp(-9\,\log n)\leq\exp(-2\log n)=\frac{1}{n^2},
\end{eqnarray*}
where  we applied Hoeffding inequality  (\ref{ho2}) with $k=2\, n^{1/2} (\log n)^{1/2}$ and $\displaystyle{x=\frac{3}{2}\left(\frac{\log n}{n}\right)^{1/4}}.$
$\Box$
\bigskip
\noindent
{\bf Proof of Theorem 1.1.}  The proof follows the basic ideas  of Theorem A.
Suppose that the walker makes $n=(cLN\log N)^2$ steps on $\bf SP(N).$  This walk can be modeled in the following way. We consider the absolute value of a simple symmetric random walk on the line $S_n.$ Then we get excursions which we throw in N urns (the legs of the spider) with equal probability. We will use Lemma 3.2 to estimate the number of tall (at least $L$ high) excursions, which are randomly placed in the $N$ urn, and then apply Theorem C. To follow this plan, let

\begin{eqnarray*}
\mu &=& \log N \\
B_n^{-}&=&\{\zeta(L,n)\leq (1-2\ep) \mu \} \\
B_n&=&\{(1-2\ep)\mu < \zeta(L,n)<(1+2 \ep) \mu \} \\
B_n^+&=&\{  \zeta(L,n)\geq (1+2\ep) \mu \}.
\end{eqnarray*}
Having

\begg
\bp(M(n,L)\}=\bp(M(n,L)|B_n^{-})\bp(B_n^{-})+\bp(M(n,L),B_n)+\bp(M(n,L)|B_n^{+})\bp(B_n^{+}),
\endd
observe that by Theorem C,
$$\lim_{N\to \infty}\bp(M(n,L)|B_n^{-})=0.$$
\newpage
\noindent
Using Lemma 3.2 we have for $n$ big enough  that
\begin{eqnarray*}
\bp(B_n)&=& \bp((1-2\ep)\mu\leq \zeta(L,n)\leq (1+2\ep)\mu) \\
&\leq& \bp \left(\frac{\xi(0,n)}{L}\leq(1+3\ep)\mu\right) - \bp \left(\frac{\xi(0,n)}{L}\leq(1-\ep)\mu\right)
+\frac{2}{n},
\end{eqnarray*}
where we used that the condition $L\log N\leq N$ of the theorem ensures that
\begg
\ep \mu\geq 4 n^{1/4} \log n^{1/4} \label{con}
\endd
 for  large enough $N.$ Consequently we have, that
$$\lim_{n\to\infty} \bp(B_n)\leq\lim_{n\to\infty} \bp\left(|Z|\leq\frac{ (1+3\ep)}{c}\right)-\bp\left(|Z|\leq\frac{ (1-\ep)}{c}\right)+\frac{2}{n}$$
Thus as $\ep \to 0$
$$\lim_{N\to \infty}\bp(B_n)=0.$$
Again by Theorem C

$$\lim_{N\to \infty}\bp(M(n,L)|B_n^{+})=1.$$
and by Lemma 3.2\, if $n$ is big enough  and $L \log N \leq N$ we have  using (\ref{con}) again that

\begin{eqnarray*}
\bp(B_n^+)&=&\bp(  \zeta(L,n)\geq (1+2\ep) \mu ) \\
&\geq&\bp\left(  \frac{\xi(0,n)}{L}\geq (1+3\ep) \mu \right)+\frac{1}{n},
\end{eqnarray*}
Consequently
$$\lim_{N \to \infty}\bp(B_n^+)\geq\bp\left(|Z|\geq \frac{1+3\ep}{c}\right).$$
and similarly

\begin{eqnarray*}
\bp(B_n^+\}&=&\bp(  \zeta(L,n)\geq (1+2\ep) \mu ) \\
&\leq& \bp\left(  \frac{\xi(0,n)}{L}\geq (1+\ep) \mu \right)+\frac{1}{n}
\end{eqnarray*}
and
$$\lim_{N \to \infty}\bp(B_n^+)\leq \bp\left(|Z|\geq \frac{1+\ep}{c}\right).$$
Letting $\ep\to 0,$ we finally get that
$$\lim_{N \to \infty}\bp(B_n^+)=\bp\left(|Z|\geq \frac{1}{c}\right)=p(c).\quad \quad\Box$$

\bigskip

\noindent
{\bf Proof of Theorem 1.2.} We use the  notations of the previous theorem, with the sole  exception that now  $n=(f(N)L N\log N)^2.$
Observe, that
\begg
\bp(M(n,L))\geq\bp(M(n,L)|B_n^{+})\bp(B_n^{+}),
\endd
and as above we know that
$$\lim_{N \to \infty}\bp(M(n,L)|B_n^{+})=1.$$
So we only have to show that
$$\lim_{N \to \infty}\bp(B_n^{+})=1.$$
Now again by Lemma 3.2
\begin{eqnarray*}
\bp(B_n^{+})&=&\bp(  \zeta(L,n)\geq (1+2\ep) \mu )\geq \bp\left(\frac{\xi(0,n)}{L}\geq 4n^{1/4}(\log n)^{3/4}+(1+2\ep) N \log N\right)- \frac{1}{n} \\
&=&\bp\left(\frac{\xi(0,n)}{n^{1/2}}\geq \frac{4n^{1/4}(\log n)^{3/4}}{N\log N f(N)}+\frac{(1+2\ep) N \log N}{N\log N f(N)}\right)- \frac{1}{n}.
\end{eqnarray*}
Being $L\log N\leq N$ and $f(N) \to \infty$, it is easy to see that

$$\lim_{N \to \infty}\frac{4n^{1/4}(\log n)^{3/4}}{N\log N f(N)}+\frac{(1+2\ep) N \log N}{N\log N f(N)}=0,$$
Thus the  limit of the above probability is  $\bp(|Z|\geq 0)=1$
which proves our theorem.
$\Box$

\bigskip
\noindent
{\bf Proof of Theorems 1.3. and 1.4.} To prove these two theorems, it is enough to  repeat the proof of Theorems 1.1 and 1.2 and apply   Lemma B instead of Lemma C.

\bigskip

\noindent
{\bf Proof of Theorem 1.5.} Notations are the same as in Theorem 1.2, except that now $f(N)\downarrow 0.$ During the proof we suppose that $f(N)L N\log N \to \infty$  otherwise there is nothing to prove. Observe that
\begg
\bp(M(n,L))\leq\bp(M(n,L)|B_n^{-})\bp(B_n^{-})+ \bp(\overline{B_n^{-}}).
\endd
As we know  from Theorem B that $\lim_{N \to \infty}\bp(M(n,L)|B_n^{-})=0,$  it is enough to prove that
$\lim_{N \to \infty} \bp(\overline{B_n^{-}})=0.$
 We show that $\lim_{N \to \infty}\bp(B_n^{-})= 1.$  Using Lemma 3.2 and the condition $L\log N\leq N,$ we have
\begin{eqnarray*} \bp (B_n^{-})&=&\bp(\zeta(L,n)\leq (1-2\ep) \mu )  \\
&\geq& \bp\left(\frac{\xi(0,n)}{L}+4n^{1/4}(\log n)^{3/4}\leq (1-2\ep)\mu \right) \\
&=& \bp\left(\frac{\xi(0,n)}{\sqrt{n}}
\leq \frac{(1-2\ep) }{ f(N)}-\frac{4n^{1/4}(\log n)^{3/4}}{N\log N f(N)}\right)\\
&\geq& \bp\left(\frac{\xi(0,n)}{\sqrt{n}}\leq \frac{1}{f(N)}\left(1-2\ep-\frac{4f^{1/2}(N) (4 \log N +2\log f(N))^{3/4}}{\log N}\right)\right).
\end{eqnarray*}
It is easy to see that  as $N\to\infty,$    $\frac{1-2\ep}{f(N)}\to +\infty,$ while the   fraction next to $(1-2\ep)$  goes to 0.
Consequently
$$\lim_{N\to\infty} \bp (B_n^{-})=P(|Z|< +\infty)=1.\quad \quad \Box$$

\bigskip
\noindent
{\bf Proof of Theorem 1.6.}  Fix the integer $K.$ The proof of this theorem will be given in three major steps. First we show that a simple symmetric walk  $S_n$ almost surely has a subsequence $S_{n_k}$ such that it has $K$ excursions which are as high as the theorem states. Then we show that $S_n$ can't have at least $K$ excursions which are all higher than it is stated. Finally we show that these imply our theorem for $\bf{SP(K)}.$
Consider the following points

$$ P_j=\left(\frac{j}{2K-1}, \frac{r_j}{2K-1}\right)\quad j=0,1,... 2K-1, \quad
{\rm where}\,\, r_j=\sin(j\frac{\pi}{2}).$$
Now define the zigzag  function $f(x)$ on $[0,1]$ by connecting the consecutive $ P_j$ points with line segments (linear interpolation).
It is easy to see that this $f(x)\in {\cal S},$ with  $\int_0^1\dot{f}^2(x)dx=1.$ Clearly $|f(x)|$ has exactly K  maximum points, and in between of any two of them $f(x)$ crosses the $x$- axis.    Applying now Theorem F
we conclude that with probability one, there exists a sequence $\{n_k=n_k(\omega)\}$ such that
$$\sup_{0\leq x \leq 1}\left\vert \frac{S(n_k x)}{\sqrt{2n_k\log\log n_k}}-f(x)\right\vert\to 0.$$
 Thus with probability one for any $\ep>0,$ and for $k$ big enough

$$\left\vert S\left(\left[n_k \frac{2j-1}{2K-1}\right]\right)\right \vert\geq (1-\ep)\frac{1}{2K-1}\sqrt{2n_k\log\log n_k}\quad  j=1,2,...   K,$$
 hence

\begg
 \limsup_{k\to\infty} \left\vert\frac{S\left([n_k \frac{2j-1}{2K-1}]\right)} {\sqrt{2n_k\log\log n_k}}\right\vert\geq\frac{1}{2K-1}
 \quad  j=1,2,...   K.
 \label{six}
\endd
and
 $$\lim_{k\to\infty} \left\vert\frac{S\left([n_k \frac{2j}{2K-1}]\right)} {\sqrt{2n_k\log\log n_k}}\right\vert=0 \quad  j=1,2,...   K-1.$$

These observations allow us to conclude that $S_{n_k}$ almost surely has $K$ excursions which are, roughly saying,
$\displaystyle{\frac{1}{2K-1} \sqrt{2n_k\log\log n_k}}\,$ tall, as the theorem claims. We  will simply call these excursions tall.

Now we  want to show that it is impossible to construct a sequence $n_k$ such that the corresponding path should contain with probability one, at least  $K$ excursions all of which are  taller than those ones  above. By Strassen theorem we know that if there would exist such a subsequence, then it would contain a further subsequence, which with the above used normalization has to converge uniformly  to a function $g(x),\,\, 0\leq x\leq 1,$  and this $g(x)$ has to be an element of ${\cal S}.$  If such a $g$ would exist within ${\cal S},$
then by Lemma G if would need to be linear, so we are looking for a linear function with at least $K$  absolute maximum points, all of which  maxima have to be bigger than $\frac{1}{2K-1}.$
Suppose that we could find such a linear function $g(x)$ which has at least K absolute maxima all of which are at least $\alpha > \frac{1}{2K-1}.$ It is obvious that for  minimizing $I(g)$ the number of such maxima should be exactly $K$ and  the value of all the maxima should be exactly $\alpha.$ To get $K$ excursions our zigzag linear function $g(x)$ always has to return to zero between consecutive maximums. All what remains to show that to minimize $I(g(x))$ all the different sections of linearity of $g(x)$ should be equally spaced as in the construction of our $f(x).$  To see that, it is enough to consider that if the length of these consecutive intervals on the $x$-axis are $x_1,  x_2, ... x_{2K-1}$  (with $x_i>0$ for all $i$) with a total length $x_1+  x_2+ ...+ x_{2K-1}=a\leq 1, $ then

$$I(g)=\alpha^2 \sum_{i=1}^{2K-1}\frac{1}{x_i}$$
So we have to solve the minimization  problem
$$\min \,\,\alpha^2 \sum_{i=1}^{2K-1}\frac{1}{x_i}$$
under the condition of $x_1+  x_2+ ...+ x_{2K-1}=a. $
It is an easy calculation to see, that the solution is
$$x_i=\frac{a}{2K-1} \quad i=1,2...2K-1$$
for which  $I(g)=\frac{\alpha (2K-1)}{a}.$
Thus to get the smallest $I(g),$  $a$ should be selected to be 1, and then selecting $\alpha=\frac{1}{2K-1} $ to have $I(g)=1$,
essentially gives back our $f(x).$  The essentially words means here, that we can select $2^{2K-1}$ zigzag functions with the given height by simply selecting which of the pikes should be on the positive side of the $x$-axis. So we proved that in (\ref{six})
we actually have equality.
\begg
 \limsup_{k\to\infty}  \left\vert \frac{S\left([n_k \frac{2j-1}{2K-1}]\right)} {\sqrt{2n_k\log\log n_k}}\right\vert=\frac{1}{2K-1}\quad  j=1,2,...   K. \label{equal}
\endd
Now returning to $\bf{SP(K)}$, we observe that to get a spiderwalk on $\bf{SP(K)}$ from  an ordinary simple symmetric walk  $S_n$, we just have to consider the consecutive excursions of $S_n$
 and put each excursion with equal probability  to one  of the $K$ legs of our spider. We have shown above the almost sure  existence of a   subsequence  $S_{n_k}, $  which has $K$ tall excursions. All what remains to  show is, that such $S_{n_k} $ has a further subsequence  $S_{n_{k_j}} $ such that when the excursions are randomly placed to the $K$ legs, all of its $K$ tall excursions are on different legs. To see this, it is enough to observe, that every time when the tall excursions of $S_{n_{k}} $  are randomly placed to the $K$ legs, then with probability $\frac{K!}{K^K}$ they all fell on different legs. Applying now Lemma G proves  the theorem. $\Box$

\end{document}